PAPER • OPEN ACCESS

# Information security control as a task of control a dynamic system

View the article online for updates and enhancements.

# Information security control as a task of control a dynamic system

S N Masaev[1,2,6], A N Minkin[1,3], Yu N Bezborodov[1,3], D A Edimichev[1] and Y K Salal[4,5]

[1] Siberian Federal University, pr. Svobodnyj, 79, Krasnoyarsk, 660041, Russia
[2] Control Systems LLC, 86 Pavlova Street, Krasnoyarsk, 660122, Russia
[3] FSBEI HE Siberian Fire and Rescue Academy EMERCOM of Russia, 1 Severnaya Street, Zheleznogorsk, 662972, Russia
[4] South Ural State University, pr. Lenina 76, Chelyabinsk, 454080, Russia
[5] University of Al-Qadisiyah, Diwaniyah, 58002, Iraq

[6] E-mail: faberi@list.ru

**Abstract**. The area of research includes control theory, dynamic systems, parameters of the external environment, mode, integral indicators, British standards. The main idea of the article is information security. The activity of a large-scale object (enterprise) is considered. The activity of the enterprise is presented as a multidimensional dynamic system and is displayed as a digital copy of 1.2 million parameters. A British digital copy-based information security standard is being introduced. Information security equipment and software were purchased. The training of the company's personnel was completed. Evaluation of implementation (activities) is done as an integral indicator. The dynamics of the integral indicator assesses the implementation of the British standard.

## 1. Introduction

Information security is the most pressing issue of the digital society and its interaction. British standard BS [1] helps to provide the necessary level of information at the enterprise in automated production control systems. Foreign authors D. Waterman, S. Russell, P. Norvig were mainly engaged in the implementation of the goals of increasing the security of automated information processes [2-4].

The common ancestor of international standards for information security is the British standard BS 7799. The UK government needed such a standard and it was developed in 1995 by his order with the name "Code of Practice for Information Security Management". This is a regulation on information security control in an organization, which includes 10 areas and 127 controls. This is a solid foundation for an ISMS. It combines the best practices from around the world. The standard turned out to be in great demand. In 1998, BS 7799-2 Specification and Application Guidance was developed. This document defines the requirements for the first part of BS 7799-1. BS 7799-2 has also received international recognition. It is used in leading enterprises in the world and countries. In 1999, this required a revision of both parts of the British standard to align with the international quality systems ISO 9001 and ISO 14001. In 2000 ISO approved BS 7799-1 as the international standard ISO / IEC 17799:2000 without modification by the technical committee.

In 2002, the second part of the BS 7799 standard was similarly revised. Three years later, ISO was adopted as the international standard ISO/IEC 27001:2005 "Information technology - Security methods - Information security management systems - Requirements" and the first part of the standard was updated.

From 1995 to 2005, the prestige of the standard has grown significantly. Creation of ISO 27001 with ISMS specifications gained international status and prominence. The increasing theft of personal information of citizens, employees and company data from information databases increases the role of information security management in accordance with the ISO 27001 standard.

Traditionally in this area there is a problem of control an object of large dimension. The large dimension of the object (enterprise) makes it difficult to effectively use even simple control approaches. Applying any control approach requires appropriate staff competencies. With the release of ISO 27001, the ISMS specifications have acquired an international status, and we can now expect a significant increase in the role and prestige of ISMS certified according to the ISO 27001 standard.

In 2006, the third part of the BS 7799 standard was developed. It is devoted to the assessment of information security risks. However, the standard does not have specific calculation methods and is not officially recognized at the international level.

V. V. Leontiev and L. V. Kantorovich, A. G. Granberg, A. G. Aganbegyan, V. F. Krotov and others were involved in the control of research objects [5-7]. The authors have formed a classical form for representing a research object (enterprise) and the methods used on them.

It is relevant to consider the impact of the acquired competencies of personnel of the BS (BS 7799-3:2006) [1] on the business processes of an object (enterprise) in an unpredictable external environment.

From a practical point of view, it is relevant to consider the introduction of a new BS information security management standard for enterprise of Russia. Since BS is more important than the ISO series standards and isn't widespread in Russia.

A purpose of research is: to estimate the state of an research object as a multidimensional dynamic system in the basic mode of operation and its control mode, through skills set of BS (BS 7799-3:2006) with unknown parameters of the external environment.

## 2. Method

It is enough to imagine the activity of the enterprise as $S=\{T,X\}$, where $T=\{t:t=1,...,T_{max}\}$ - a lot of time points with a selected interval for analysis; $X$ - set of system parameters; $x(t)=[x^1(t),x^2(t),...,x^n(t)]^T \in X - n$ – vector of indicators corresponding to the state of the system. Indicators of the vector $x^i(t)$ - the value of financial expenses and income of the enterprise. The dimension of the system $n$ is 1.2 million parameters. Based on the parameters $X$ and $T$, we consider our research object a multidimensional dynamic system (hereinafter referred to as the system). The system has different dimensions at each moment of time and is a digital copy of the enterprise.

The analysis of the system at the moment $t$ is performed $x(t)$ on the basis $k$ of previous measures. The parameter $k$ is the length of the time series segment (accepted $k=6$ months in the work). Then we have a matrix

$$X_k(t) = \begin{bmatrix} x^T(t-1) \\ x^T(t-2) \\ \dots \\ x^T(t-k) \end{bmatrix} = \begin{bmatrix} x^1(t-k) & x^2(t-k) & \cdots & x^n(t-k) \\ x^1(t-k) & x^2(t-k) & \cdots & x^n(t-k) \\ \dots & \dots & \dots & \dots \\ x^1(t-k) & x^2(t-k) & \cdots & x^n(t-k) \end{bmatrix} \quad (1)$$

$$R_k(t) = \frac{1}{k-1} \overset{o}{X}_k^T(t) \overset{o}{X}_k(t) = \left\| r_{ij}(t) \right\|, \quad (2)$$

$$r_{ij}(t) = \frac{1}{k-1} \sum_{l=1}^{k} \overset{o}{x^i}(t-l) \overset{o}{x^j}(t-l), \quad i,j = 1,\ldots,n, \tag{3}$$

where $t$ are the time instants, $r_{ij}(t)$ are the correlation coefficients of the variables $x^i(t)$ и $x^j(t)$ at the time instant $t$.

Next we form one of the four integral indicators – the sum of the absolute indicators of the correlation coefficients. It is indicator for express estimation of the correlation of system parameters $G_i(t)$:

$$R_i(t) = G_i(t) = \sum_{j=1}^{n} |r_{ij}(t)|. \tag{4}$$

The state of the entire system is calculated as:

$$G = \sum_{t=1}^{T=\max} \sum_{i=1}^{n} G_i(t). \tag{5}$$

The $V$ is the set of BS (BS 7799-3:2006), which can be represented $V_i^k$ as a set of employee skills of BS:

$$V = \sum_{t=1}^{T=\max} \sum_{i=1}^{n} V_i^k(t). \tag{6}$$

We will carry out the identification of the performed functions of the system with the set of employee skills of BS (BS 7799-3:2006) for the function (business-process). Each skill of BS $V_i^k$ is characterized by the business-process of the enterprise:

$$V_i^k = \sum_{i=1}^{n} v_i^j(x_j^i) \to \min, \tag{7}$$

where $v_i^j$ is name of a skill (compliance $x_j^i$ is $v_i^j$ set as 1-yes, 0-no); $x_j^i$ - the costs of the $i$ - the skill of BS (BS 7799-3:2006) for the function (business-process) $j$. Control method is $V$. It is skill of employee $v_i^j$ so have a vector of strategic planning $v(t) = [v^1(t), v^2(t), \ldots, v^n(t)]^T \in V$ - $n$ - dimension. Then $V = \sum_{t=1}^{T=\max} \sum_{i=1}^{n} V_i^k(t)$.

Payment of the functional duties of employees of the system is limited by resources $C$, then $C(X) \leq C$. This restriction applies to all subsystems of the researched system.

The implementation of the method is performed in the author's complex of programs [8,9].

## 3. Characteristics of the research objects

The considered economic object (enterprise) belongs to the woodworking industry of the of the Krasnoyarsk Territory. The volume of roundwood processing is 800 thousand cubic meters. The company employs 600 people. The company has a standard management structure through the organizational system: economic department, production department, logistics and transport department, finance department, accounting department, sales department, marketing department. In 1.5 years, it is planned to double production at the expense of bank loans, taking into account tax benefits under certain scenarios of market development and enterprise strategy.

At this enterprise, a management loop is being implemented through the competencies of personnel in compliance with BS 7799-3:2006 [1]: the organization must assess whether it has enough resources to perform the following actions: perform a risk assessment and develop an information base (IB) risk treatment plan, determine and implement policies and procedures, including the implementation of the selected controls, monitor information security risk management tools, monitor the information security risk management process. etc.

### 4. Experiment result

Initial data in the calculation algorithm: $X=5{,}641{,}442$ thousand rubles, $n=1.2$ million values, the control loop is $V_{slills\_BS}$ - skills set of BS (BS 7799-3:2006). The calculation of the standard operating mode of the enterprise $V_{(basic\ mode)}=5{,}069.93$,

From the 1st period, the personnel in the departments are trained according to skills set of BS. Additional servers are purchased to control information flows. The software is purchased for organizing access rights to information flows between business processes of an enterprise. Personnel are trained to control all important points of aggregation of production information. Also, a business trip for training personnel in skills set of BS (BS 7799-3:2006) of the enterprise (dynamic system) is paid.

The simulation algorithm execution time is 436 minutes [8,9].

A table 1 shows the experiment result of estimating the control mode $V_i(t)$ through skills set of BS (BS 7799-3:2006).

**Table 1.** Regimes: $V_{(basic\ mode)}$ and $V_{(skills\_BS)}$.

| $t$ | $V_{(basic\ mode)}$ | $V_{(skills\_BS)}$ | $\Delta V$ | $t$ | $V_{(basic\ mode)}$ | $V_{(skills\_BS)}$ | $\Delta V$ |
|---|---|---|---|---|---|---|---|
| 1 | 87.34 | 110.57 | 23.23 | 30 | 96.32 | 95.32 | -1.00 |
| 2 | 70.94 | 90.15 | 19.21 | 31 | 105.10 | 104.10 | -1.00 |
| 3 | 51.43 | 69.07 | 17.64 | 32 | 98.66 | 97.66 | -1.00 |
| 4 | 56.35 | 81.21 | 24.86 | 33 | 82.19 | 81.19 | -1.00 |
| 5 | 59.26 | 89.64 | 30.37 | 34 | 76.23 | 75.22 | -1.00 |
| 6 | 73.39 | 110.67 | 37.29 | 35 | 68.52 | 68.52 | 0.00 |
| 7 | 95.25 | 123.53 | 28.29 | 36 | 60.51 | 60.51 | 0.00 |
| 8 | 92.64 | 124.07 | 31.43 | 37 | 53.13 | 53.13 | 0.00 |
| 9 | 95.53 | 125.43 | 29.90 | 38 | 61.65 | 61.65 | 0.00 |
| 10 | 70.17 | 100.32 | 30.15 | 39 | 53.51 | 53.51 | 0.00 |
| 11 | 58.42 | 74.76 | 16.35 | 40 | 51.84 | 51.84 | 0.00 |
| 12 | 56.48 | 72.97 | 16.50 | 41 | 72.03 | 72.03 | 0.00 |
| 13 | 61.88 | 79.69 | 17.81 | 42 | 93.08 | 93.08 | 0.00 |
| 14 | 71.87 | 90.79 | 18.92 | 43 | 99.23 | 99.23 | 0.00 |
| 15 | 52.45 | 67.93 | 15.48 | 44 | 115.79 | 115.79 | 0.00 |
| 16 | 53.92 | 70.08 | 16.16 | 45 | 110.12 | 110.12 | 0.00 |
| 17 | 84.06 | 101.54 | 17.48 | 46 | 103.64 | 103.64 | 0.00 |
| 18 | 114.43 | 132.62 | 18.19 | 47 | 88.22 | 87.22 | -0.99 |
| 19 | 132.20 | 132.20 | 0.00 | 48 | 69.27 | 68.27 | -1.00 |
| 20 | 153.90 | 153.92 | 0.01 | 49 | 55.14 | 54.14 | -1.00 |
| 21 | 164.54 | 164.53 | -0.01 | 50 | 63.51 | 62.51 | -1.00 |
| 22 | 150.02 | 151.00 | 0.98 | 51 | 50.02 | 49.02 | -1.00 |
| 23 | 140.74 | 144.74 | 4.00 | 52 | 61.58 | 60.58 | -1.00 |
| 24 | 115.09 | 119.09 | 3.99 | 53 | 60.33 | 59.33 | -1.00 |
| 25 | 87.02 | 91.02 | 4.00 | 54 | 147.33 | 147.33 | 0.00 |
| 26 | 100.60 | 104.60 | 4.01 | 55 | 158.41 | 158.41 | 0.00 |
| 27 | 87.59 | 91.59 | 4.00 | 56 | 156.87 | 156.87 | 0.00 |
| 28 | 76.04 | 79.04 | 3.00 | 57 | 167.90 | 167.90 | 0.00 |
| 29 | 76.26 | 76.26 | 0.00 | **Total** | 5,069.93 | 5,491.17 | 421.24 |

A figure 1 shows the experiment result of estimating the control mode $V_i(t)$ through the skills set of BS (BS 7799-3:2006).

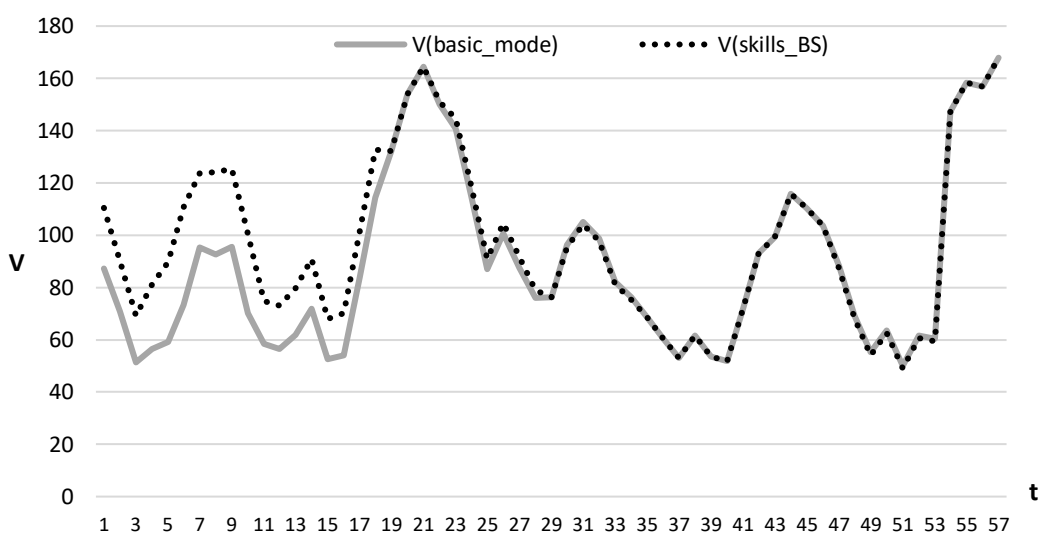

**Figure 1.** Indicator dynamics $V_i(t)$.

## 5. The discussion of the results
BS 7799-3: 2006 has important characteristics of flexibility and versatility. The approaches formulated in it are generally applicable to all organizations. This is not influenced by either the form of ownership, nor the geographical position, nor the form of detail. The standard regulates any technology used in the enterprise.

The standard defines a list of issues for development: the beginning of the implementation of an information security system, how to regulate it, audit criteria. It is very well suited to justify an enterprise's information security costs. The standard is the authoritative source for this.

The standard requires revision at each implementation enterprise, since only basic information security issues are written in it. Insufficient adaptation of the standard to the specifics of the enterprise becomes a threat of information loss.

The enterprise must determine the line between the detail of the implementation of the standard and the measure of information security for itself.

## 6. Conclusion
The BS 7799-3:2006 does not limit the depth of security controls and does not impose technical implementation. However, there is a reference to specific technologies: USB keys, smart cards, certificates, etc. It is important to understand that the standard does not compare technologies with each other and does not highlight the best one.

Nevertheless, the research was able to assess the implementation of BS 7799-3:2006 at the enterprise with the exact value of the integral indicator $V_{(skills\_BS)}$ – 5,491.17. The dynamics of the integral indicator displays all the changes during the implementation of BS 7799-3:2006: wages, business trips, taxes, personnel training, office, communications, org. technics. It costs an additional 9,060 thousand rubles. The total costs of the enterprise for five years with the implementation of skills set of BS (BS 7799-3:2006) will amount to 5,666,745 thousand rubles.

Therefore, the assessment of the transition to object control through skills set of BS (BS 7799-3:2006) is estimated as $\Delta V = V_{(skills\_BS)} - V_{(basic\_mode)} = 421.24$. The purpose of the research has been achieved.